\documentclass{amsart}
\usepackage{}
\usepackage{color}
\usepackage{amssymb}
\usepackage{amsthm}

\newcommand{\ul}{\underline}

\newtheorem{theorem}{Theorem}[section]
\newtheorem{lemma}[theorem]{Lemma}

 \theoremstyle{definition}

\theoremstyle{remark}
\newtheorem{remark}[theorem]{Remark}

\numberwithin{equation}{section}

\begin{document}

\title[The Dirichlet problem]
{The Dirichlet problem for Monge-Amp\`ere type equations on Riemannian manifolds}

\author{Weisong Dong}
\address{School of Mathematics, Tianjin University, 135 Yaguan Road,
Jinnan, Tianjin, China, 300354}
\email{dr.dong@tju.edu.cn}

\author{Jinling Niu}
\address{Department of Mathematics and Statistics, Changshu Institute of Technology, Changshu, China, 215500}

\email{njlhit@163.com}

\author{Nadilamu Nizhamuding}
\address{School of Mathematics, Tianjin University, 135 Yaguan Road,
Jinnan, Tianjin, China, 300354}
\email{adira@tju.edu.cn}

%\date{}

\begin{abstract}

In this paper, we study the Dirichlet problem for Monge-Amp\`ere type equations for $p$-plurisubharmonic functions on Riemannian manifolds. The $a$   $priori$  estimates up to the second order derivatives of solutions are established. The existence of a solution then follows by the continuity method.

\emph{Mathematical Subject Classification (2020):} 35J15; 35B45.

\emph{Keywords:}  Monge-Amp\`ere  type equations;  $a$  $priori$  estimates; Riemannian manifolds; existence.

\end{abstract}

\maketitle

\section{Introduction}

Let $(M, g)$ be a compact Riemannian manifold of dimension $n \geqslant 2$ with smooth boundary $\partial M$ and $\bar M := M \cup \partial M$. Let $\nabla^2 u$ denote the Hessian of a function $u$ and $\chi$ be a $(0, 2)$ tensor on $M$. Given any smooth function $u\in C^\infty (M)$, the operator  $\mathcal{M}_p^n (u)$ is defined as below
$$\mathcal{M}^n_p(u)=\Pi_{1\leqslant i_1 < \cdots < i_p\leqslant n} (\lambda_{i_1} + \cdots + \lambda_{i_p}),$$
where $\lambda=(\lambda_1, \cdots, \lambda_n)$ are the eigenvalues of $U:= \nabla^2 u + \chi$ with respect to the metric $g$. In this paper we are concerned with the Dirichlet problem for  Monge-Amp\`ere type equation of the form
\begin{equation}
\label{eqn}
\mathcal{M}^n_p(u) = f (x, u, \nabla u)\;  \mbox{in}\; M,
\end{equation}
with the Dirichlet  boundary data
\begin{equation}
\label{eqn-b}
u=\varphi \; \mbox{on}\; \partial M,
\end{equation}
where $f\in C^\infty(M \times \mathbb{R} \times \mathbb{R}^n)>0$  and $\varphi \in C^\infty(\partial M )$  is a given function.
The Monge-Amp\`ere equation, that corresponds to  $p=1$ in \eqref{eqn}, was studied in Euclidean space by Caffarelli-Nirenberg-Spruck \cite{CNS1}, which played a pivotal role in many geometry problems such as the Minkowski problem (\cite{CY}, \cite{N}, \cite{P}) and the Weyl problem (\cite{GLL1}, \cite{HZ}, \cite{N}). Atallah-Zuily \cite{AZ} solved the Monge-Amp\`ere equation in Euclidean space with arbitrary Riemannian metric for $f=f(x)$. And in the case where  $f=f(x,u, \nabla u)$, it  was studied by Atallah  \cite{A} assuming the existence of upper and subsolutions. The Poisson equation corresponding to  $p=n$ in \eqref{eqn}, was studied by Gilbarg-Trudinger \cite{GT} and Ladyzhenskaya-Uraltseva \cite{LU}.

The equation \eqref{eqn} with $p < n$ is a fully nonlinear elliptic equation for $p$-plurisubharmonic functions.
A counterpart to \eqref{eqn} that was explored extensively is the $k$-Hessian equation as follows
\begin{equation}
\label{eqn-c}
\sum_{1\leqslant i_1 < \cdots < i_k\leqslant n} \lambda_{i_1} \cdots \lambda_{i_k}=f,
\end{equation}
where $1\leqslant k\leqslant n$. It is easy to see that when $k=n$  and $k=1$,  \eqref{eqn-c} are the Monge-Amp\`ere equation and the Poisson equation, respectively.
When $f$ depends on $\nabla u$, Guan-Ren-Wang \cite{GRW} proved the second order estimates for convex solutions. Ren-Wang \cite{RW1,RW2} removed the convexity condition for $k=n-1$ and $k=n-2$.
Li in \cite{YYLi} initiated the study of the equation \eqref{eqn-c} on closed Riemannian manifolds. He also found a necessary and sufficient condition for the existence of a solution to the Monge-Amp\`ere equation on the flat torus $\mathbb{T}^n$. Recently, Guo-Song established a necessary and sufficient condition for solving general Hessian
equations on closed real or complex manifolds in \cite{GS}. The quotient type of equation \eqref{eqn-c} arises from complex geometry which is called the $J$-equation has been studied by many authors, see \cite{Chen, DP, SW} and references therein.

 For the Dirichlet problem, B. Guan and his collaborators have done many excellent works on very general fully nonlinear elliptic equations including \eqref{eqn} and \eqref{eqn-c}. See \cite{G23} for a more recent work. Especially, Guan proved the second order estimates for \eqref{eqn} on Riemannian manifolds in \cite{GUAN2,GJ} where $f$ is assumed to be convex with respect to $\nabla u$. The same condition is also required in the work of Jiang-Trudinger-Yang \cite{JTY} to solve the Dirichlet problem for a class of augmented Hessian equations. Recently, Jiao-Liu  \cite{JL1} treated the Dirichlet problem of \eqref{eqn} with $p=n-1$ without the convex condition using the idea of \cite{CJ}. In complex settings, the operator $\mathcal{M}_{n-1}$ is related to the Gauduchon conjecture which was resolved by Sz\'{e}kelyhidi-Tosatti-Weinkove \cite{S-T-W}. The operator $\mathcal{M}_{n-1}$ also appears in the study of the form type Calabi-Yau equations raised by Fu-Wang-Wu \cite{FWW1}, which was studied subsequently by Fu-Wang-Wu \cite{FWW2} and Tosatti-Weinkove \cite{TW1}. In \cite{D2}, the first named author proved the Dirichlet problem of \eqref{eqn} with $p=n-1$ for the case $f=f(x)$ on Hermitian manifolds by deriving a second order estimate with the form as in Hou-Ma-Wu \cite{HMW2} following the argument in Collins-Picard \cite{CP} where they proved the Dirichlet problem for \eqref{eqn-c} on Hermitian manifolds.
 Here, based on many properties of the operators $\mathcal{M}^n_p(u)$ proved in \cite{Dinew,Dong,Dong1}, we investigate the Dirichlet problem for the equation \eqref{eqn} with $p\geqslant \frac{n}{2}$ on Riemannian manifolds without the assumptions about the dependence of  $f$ on $\nabla u$.

Now, we define admissible functions, for which the equation \eqref{eqn} is elliptic.
Recall that the $\mathcal{P}_p$ cone for $p = 1, \cdots, n$ is defined by
$$ \mathcal{P}_p = \{\lambda\in \mathbb{R}^n: \;  \forall \; 1\leqslant i_1 < \cdots < i_p\leqslant n, \; \lambda_{i_1} + \cdots + \lambda_{i_p} > 0\}.$$
The cone $P_p$ of symmetric $n\times n$ matrices is defined as $A\in P_p$ if $\lambda(A) \in \mathcal{P}_p$.
A function  $u \in C^2 (M)$ is called admissible (or $p$-plurisubharmonic) if $\lambda(\nabla ^2 u+\chi)(x)\in\mathcal{P}_p$ for all $x \in M$.
We refer the reader to \cite{HL13} for $p$-convex cones and to \cite{Wu,Sha1} for more results about $p$-convexity in geometry.
It is also worth to remark that constant rank theorems for $p$-convex solutions to some semilinear equations were obtained by Han-Ma-Wu \cite{HMW}.
In this paper, we suppose that there exists an admissible subsolution $\ul u \in  C^\infty(\bar M)$  satisfying
\begin{equation}
\label{26}
\mathcal{M}^n_{p}(\ul u) \geqslant f (x,\ul u, \nabla \ul u) \; \mbox{in} \; M
\end{equation}
and
$\ul u = \varphi$ on $\partial M$. The subsolution is crucial when deriving $|\nabla u|_{C^0(\partial M)}$
and $|\nabla^2 u|_{C^0(M)}$.
We also need the following growth condition
\begin{equation}
\label{27}
| \tilde f_{x_i}| + | \tilde f_z| |\nabla u| + | \tilde f_{p_l}| |\nabla u|^2 \leqslant  C|\nabla u|^{2+\gamma_0},
\end{equation}
where $\tilde f = f^{ \frac{1}{C_n^p} }$,  $0< \gamma_0 <1$ and $C$ is a uniform constant. This condition will be used to derive the gradient estimate. Actually, the gradient estimate can be proved similarly by the argument in \cite{GUAN2} since the operator $\mathcal{M}^n_p (u)$ satisfies the fundamental structure conditions therein. Combining the $a$ $priori$ estimates, in particular, Theorem 3.3, Theorem 4.1 and Theorem 5.1 and the continuity method in \cite{GUAN2}, we can prove the existence result as follows.
\begin{theorem}
Let  $f \in C^\infty(M\times \mathbb{R} \times \mathbb{R}^n)$  be a positive function  satisfying  \eqref{27}  with  $0< \gamma_0 <1$ and
\begin{equation}
\label{36}
\sup_{(x,z,p) \in  T^{\ast} M \times \mathbb{R} } \frac{- f_z(x,z,p)}{ f(x,z,p)} < \infty.
\end{equation}
Assume that $p \geqslant \frac{n}{2}$ and there is an admissible subsolution $\ul u\in C^\infty (\bar M)$ satisfying \eqref{26} and \eqref{eqn-b}.
Then, there exists an admissible solution $ u\in C^\infty (\bar M)$  of \eqref{eqn} and  \eqref{eqn-b}. Moreover, the solution is unique if $f_z\geqslant 0.$
\end{theorem}

\begin{remark}
It remains an interesting question to prove the above theorem for  $1 < p < \frac{2}{n}$.
\end{remark}

As long as we prove $|u|_{C^2 (M)} \leqslant C$,
by \eqref{36} and the degree theory in Guan \cite{GUAN2}, we immediately obtain the theorem, since $C^{2, \alpha}$ estimates follow from the Evans-Krylov theorem and more higher estimates follow from Schauder theory. By the assumption of the existence of a subsolution, $C^0$ and $C^1$ estimates are easy to derive. Note that there are no assumptions about the dependence of the right hand side function $f$ on $\nabla u$ which creates substantial difficulties to
derive the global $C^2$ estimate, as the bad term $-CU_{11}$ appears when one applies the maximum principle to the test function $W$ defined in Section 4. This also causes trouble for us to deal with third order terms. To overcome these difficulties, one can impose various assumptions on $f$ as in \cite{I1,GG}. Due to the special nature of the operator of equation \eqref{eqn}, we use idea from \cite{CJ,Dong} to control $-CU_{11}$ by good terms. Then the bad third order terms can be eliminated by the terms $-F^{1i,i1}$. One can see the  Lemma 2.4 and Lemma 2.5 in \cite{Dong1}.

The rest of the paper is organized as follows. In Section 2, we introduce some useful notations. In section 3, we derive the interior gradient estimate for a class of more generalized equation including \eqref{eqn} and present the global gradient estimate for \eqref{eqn}.  Finally, the global and boundary $C^2$ estimates are proved in Section 4 and Section 5, respectively.

\section{Preliminaries}

Throughout this paper, $\nabla$ denotes the Levi-Civita connection of $(M,g)$ and $R$ is the curvature tensor defined by
$$ R(X,Y) Z = -\nabla_X \nabla_Y Z + \nabla_Y \nabla_X Z + \nabla_{[X,Y]} Z.$$
Let $e_1, \cdots, e_n$  be a local frame on $M$ and denote $g_{ij} = g(e_i,e_j)$, $\{g^{ij}\} = \{g_{ij}\}^{-1}$, while the Christoffel symbols $\Gamma_{ij}^k$ and curvature coefficients are defined respectively by
$$\nabla_{e_i}e_j = \Gamma_{ij}^k e_k, \\\ R_{ijkl} = \langle R(e_k,e_l)e_j, e_i \rangle, \\\  R_{jkl}^i = g^{im}R_{mjkl}.$$
We shall write $\nabla_i = \nabla_{e_i}, \nabla_{ij} = \nabla_i\nabla_j - \Gamma_{ij}^k\nabla_k$, etc.

For a differentiable function $v$ defined on $M$, we identify $\nabla v$ with its gradient, and $\nabla^2 v$ denotes its Hessian which is locally given by $\nabla_{ij} v = \nabla_{i}(\nabla_{j}v) - \Gamma_{ij}^k \nabla_{k}v$. We recall that $\nabla_{ij} v = \nabla_{ji} v$,
\begin{equation}
\begin{aligned}
\label{18}
\nabla_{ijk} v-\nabla_{jik} v= R_{kij}^l \nabla_{l}v,
\end{aligned}
\end{equation}
and
\begin{equation}
\begin{aligned}
\label{25}
\nabla_{ijkl} v-\nabla_{klij} v=  &\;R_{ljk}^{m} \nabla_{im}v + \nabla_{i} R_{ljk}^{m} \nabla_{m}v +  R_{ljk}^{m} \nabla_{jm}v \\
                                  &\; + R_{jik}^{m} \nabla_{lm}v + R_{jil}^{m} \nabla_{km}v + \nabla_{k}  R_{jil}^{m} \nabla_{m}v.
\end{aligned}
\end{equation}
Define $F(\lambda)$ on the cone $\mathcal{P}_p$
by
\[
  F(\lambda) = (\Pi_{1\leqslant i_1 < \cdots < i_p\leqslant n} (\lambda_{i_1} + \cdots + \lambda_{i_p}) )^{\frac{1}{C_n^p}},
\]
where $C_n^p = \frac{n!}{p! (n-p)!}$.

Equation \eqref{eqn} then can be rewritten as
\begin{equation}
\label{eqn'}
 F (U) := F(\lambda(U)) = \tilde f(x,u, \nabla u),
\end{equation}
where  $\tilde f = f^{ \frac{1}{C_n^p} }$. We remark that $F(\lambda)$ satisfies the following structure conditions:
\begin{equation}
\label{28}
 F>0 \; \mbox{in} \; \mathcal{P}_p \; \mbox{and}  \;F=0 \;\mbox{on} \; \partial \mathcal{P}_p;
\end{equation}
\begin{equation}
\label{29}
 F_i\equiv \frac{\partial F}{\partial \lambda_i}>0 \;  \mbox{in}  \;\mathcal{P}_p,  \; \forall \; 1\leqslant i \leqslant  n;
\end{equation}
\begin{equation}
\label{30}
 F\; \mbox{is concave in} \; \mathcal{P}_p;
\end{equation}
\begin{equation}
\label{31}
F\; \mbox{is homogeneous of degree one, i.e.} \; F(t\lambda) = tF(\lambda), \; \forall \; t>0;
\end{equation}
\begin{equation}
\label{45}
\sum F_i (\lambda) \lambda_i > 0, \; \forall \; \lambda \in \mathcal{P}_p.
\end{equation}
We refer the reader to \cite{Dinew} for the properties of the operator, or to the appendix in \cite{Dong}. By  \eqref{30} and \eqref{31}, we can derive that
\begin{equation}
\label{32}
 \sum F_i(\lambda) = F(\lambda) + \sum  F_i(\lambda) (1-\lambda_i) \geqslant  F(1,\ldots,1),
\end{equation}
where $\textbf{1} = (1,\ldots,1) \in \mathcal{P}_p$.
%Without loss of generality, we may assume that $F$ is normalized such that %$F(\textbf{1}) = 1$.

Under a local frame $e_1, \cdots, e_n$, let
\[
U_{ij} = (\nabla ^2 u + \chi) (e_i,e_j) = \nabla_{ij} u + \chi_{ij}
\]
and denote
\[
F^{ij} = \frac{\partial F}{\partial U_{ij}}(U), \; \; F^{ij,kl} = \frac{\partial^2 F}{\partial U_{ij} \partial U_{kl}}(U).
\]
The matrix $\{F^{ij}\}$ has eigenvalues $F_1, \ldots, F_n$ and is positive definite by assumption \eqref{29}. And \eqref{30} implies that $F$ is concave  with respect to $U_{ij}$ by Lemma 1.13 and Corollary 1.14 in \cite{Dinew}. The equation is elliptic as the matrix $\{\frac{\partial F}{\partial U_{ij}}\}$ is positive definite for $\{U_{ij}\} \in P_p$. Moreover, when $\{U_{ij}\}$ is diagonal, so is $\{ F^{ij} \}$ and the following identities hold
\[
F^{ij} U_{ij} = \sum F_i \lambda_i, \;  \; F^{ij} U_{ik}U_{kj} = \sum F_i \lambda_i^2,
\]
where $\lambda(U) = (\lambda_1,\ldots,\lambda_n)$.

\section{A priori gradient estimates}

In this section, we establish the interior gradient estimate for a class of more generalized equation as follows
\begin{equation}
\label{eqn-2}
 F (\nabla ^2 u + A[u]) := F (\lambda (\nabla ^2 u + A[u]) )= \tilde f(x, u, \nabla u),
\end{equation}
where $(0, 2)$ tensor $A[u] = A(x, u, \nabla u)$ which may depend on $u$ and $\nabla u$. The method we applied in this section is mainly from \cite{D3}. Together with \eqref{27} we also need the following growth conditions
\begin{equation}
\label{46}
|p|^{-1} |\nabla_x A^{\xi \eta} (x, z, p)| + |A_z^{\xi \eta} (x, z, p)| \leqslant \bar {\omega}_1 (x,z) |\xi| |\eta| (1+|p|^{\gamma}), \;\;  \forall \; \xi, \eta \in T_x \bar{M}
\end{equation}
and
\begin{equation}
\begin{aligned}
\label{47}
&\;|\tilde f| + | p \cdot D_p \widetilde f(x, z, p)| + |A^{\xi\eta} (x, z, p)| / |\xi| |\eta|\\
 + &\; |p \cdot D_p A^{\xi\eta}(x, z, p)| / |\xi| |\eta| \leqslant  \bar {\omega}_2 (x,z) (1 + |p|^{\gamma}), \;\;  \forall \;\xi, \eta \in T_x \bar{M}
\end{aligned}
\end{equation}
for some continuous function $\bar{\omega}_1, \bar{\omega}_2 \geqslant 0$ and constant $\gamma \in (0,2)$ when $|p|$ is sufficiently large. Our main result is stated as follows.
\begin{theorem}
Assume that \eqref{27}, \eqref{46}, \eqref{47} hold for $0 < \gamma_0 < 1, 0 < \gamma < 2 $. Let $u \in C^3 (M)$ be an admissible solution of \eqref{eqn-2}. For any geodesic ball $B_r(x_0) \subset M$ with center to $x_0 \in M$ and radius $r < 1$,  there exists a positive constant $C$ only depending on $|u|_{C^0}$ and other known data such that
$$
| \nabla u(x_0)| \leqslant \frac{C}{r}.
$$
\end{theorem}
\begin{proof}
Without loss of generality, we suppose $0 \in M$. Choose a positive constant $r < 1$ such that $B_r(0) \subset M$. Consider the following  function on $B_r(0)$,
$$
G(x) = |\nabla u| h(u) \zeta(x),
$$
where $\zeta(x) = r^2 - \rho^2(x)$ and $\rho(x)$ denotes the geodesic distance from $0$. We can assume $r$ to be small enough so that $\rho$ is smooth and $|\nabla \rho| \equiv 1$ in $B_r(0)$. Note that $\nabla_{ij} \rho^2 (0) = 2 \delta_{ij}$. For sufficiently small $r$, we may further assume in $B_r(0)$ that
$$
\delta_{ij} \leqslant \nabla_{ij} \rho^2 \leqslant 3\delta_{ij}.
$$
Then we assume that $G(x)$ attains its maximum at an interior point $x_0 \in B_r(0)$. Choose a smooth orthonormal local frame $e_1, \ldots, e_n$ about $x_0$ such that
$$\nabla_1 u(x_0) = |\nabla u(x_0)|.$$
We may also assume $\nabla_i e_j = 0$ at $x_0$ for $i, j = 1, \ldots, n$. We differentiate log $G$ at $x_0$ twice to get that
\begin{equation}
\label{37}
\frac{\nabla_{i1} u} {\nabla_{1} u} + \frac{\nabla_{i} h} {h} + \frac{\nabla_{i} \zeta} {\zeta} = 0
\end{equation}
and
\begin{equation}
\begin{aligned}
\label{38}
0 \geqslant &\; \frac{\nabla_{ij1} u} {\nabla_{1} u} + \frac{\nabla_{ik} u \nabla_{jk} u} {(\nabla_{1} u)^2} - \frac{2\nabla_{i1} u \nabla_{j1} u} {(\nabla_{1} u)^2}   \\
&\; + \frac{\nabla_{ij} h} {h} - \frac{\nabla_{i} h \nabla_{j} h} {h^2} +  \frac{\nabla_{ij} \zeta} {\zeta} - \frac{\nabla_{i} \zeta \nabla_{j} \zeta} {\zeta^2}.
\end{aligned}
\end{equation}
Differentiating the equation \eqref{eqn} we have, at $x_0$,
$$
F^{ij} \nabla_{1ij} u + F^{ij} \nabla_{1}^{'} A^{ij} + F^{ij} A_z^{ij} \nabla_{1} u + F^{ij} A_{pl}^{ij} \nabla_{1l} u = \widetilde f_{x_1} + \widetilde f_{z} \nabla_{1} u + \widetilde f_{pl} \nabla_{1l} u.
$$
By \eqref{18} and the above equation, we obtain
\begin{equation}
\label{39}
\begin{aligned}
F^{ij} \nabla_{ij1} u = &\; \widetilde f_{x_1} + \widetilde f_{z} \nabla_{1} u + \widetilde f_{pl} \nabla_{1l} u + F^{ij} R^l_{ji1} \nabla_{l} u \\
& \; - F^{ij} \nabla_{1}^{'} A^{ij} - F^{ij} A_z^{ij} \nabla_{1} u - F^{ij} A_{pl}^{ij} \nabla_{1l} u.
\end{aligned}
\end{equation}
Combining  \eqref{27}, \eqref{46}, \eqref{47} and \eqref{39} we have
\begin{equation}
\label{40}
\begin{aligned}
\frac {F^{ij} \nabla_{ij1} u } {\nabla_{1} u}
%\geqslant &\; -C |\nabla_{1} u|^{1 + \gamma_0} (1 + \sum  F^{ii}) + %\widetilde f_{pl} \frac{\nabla_{1l} u} {\nabla_{1} u} \\
%&\; -C |\nabla_{1} u|^\gamma \sum F^{ii} - F^{ij} A_{pl}^{ij} %\frac{\nabla_{1l}} {\nabla_{1} u}\\
\geqslant &\; -C |\nabla_{1} u| ^{1+\gamma_0} (1 + \sum  F^{ii})- \frac {h^{'}} {h} \widetilde f_{p1} \nabla_{1} u + \frac{h^{'}} {h} F^{ij} A_{p1}^{ij} \nabla_{1} u \\
&\; -\frac{2Cr} {\zeta} \nabla_{1} u - C |\nabla_{1} u|^\gamma \sum  F^{ii} - \frac{2Cr} {\zeta} |\nabla_{1} u|^{\gamma-1} \sum  F^{ii},\\
\end{aligned}
\end{equation}
where we used the fact  $| \nabla \rho | \equiv 1$  and $| \widetilde f_{pl} | \leqslant C |\nabla u| $ when $|\nabla u|$ is sufficiently large, and $C$ is a uniform positive constant. From equality \eqref{37} we have
\begin{equation}
\label{41}
\begin{aligned}
F^{ij} \frac{\nabla_{i1} u \nabla_{j1} u} {(\nabla_{1} u)^2} = &\; F^{ij} \left( \frac{\nabla_{i} h}{h} + \frac{\nabla_{i} \zeta } {\zeta} \right) \left( \frac{\nabla_{j} h} {h} + \frac{\nabla_{j} \zeta } {\zeta} \right ) \\
\leqslant &\; 2 F^{ij} \frac{\nabla_{i} h} {h} \frac{\nabla_{j} h} {h} + 2 F^{ij} \frac{\nabla_{i} \zeta} {\zeta} \frac{\nabla_{i} \zeta} {\zeta}.
\end{aligned}
\end{equation}
It is easy to see that
\begin{equation}
\label{42}
\begin{aligned}
F^{ij} \frac{\nabla_{ij} \zeta} {\zeta} - F^{ij} \frac{\nabla_{i} \zeta \nabla_{j} \zeta} {\zeta^2} \geqslant  -\frac { 3 \sum F^{ii} }
{\zeta} - \frac{4 \rho^2 F^{ij} \nabla_i \rho \nabla_j \rho} {\zeta^2}.
\end{aligned}
\end{equation}
Since $\{ F^{ij} \} > 0$, contracting \eqref{38} with $\{ F^{ij} \}$ and by \eqref{37}, \eqref{40}, \eqref{41} and \eqref{42}, we get
\begin{equation}
\label{43}
\begin{aligned}
0 \geqslant  & \; -C |\nabla_{1} u|^{1 + \gamma_0} ( 1 + \sum F^{ii} ) - \frac{h^{'}} {h} \widetilde f_{p1} \nabla_{1} u + \frac{h^{'}}{h} F^{ij} A_{p1}^{ij} \nabla_{1} u\\
& \; -\frac{2Cr} {\zeta} \nabla_{1} u - \frac{2Cr} {\zeta} |\nabla_{1} u|^{\gamma - 1} \sum F^{ii} + \frac{h^{'}} {h} F^{ij} \nabla_{ij}u\\
 & \; + \left( \frac{h^{''}} {h} - 3 \frac{h^{'2}} {h^2} \right) F^{11} (\nabla_{1} u)^2 - \frac{3} {\zeta} \sum F^{ii}
 - \frac{12r^2} {\zeta^2} F^{ij} \nabla_{i} \rho \nabla_{j} \rho,
\end{aligned}
\end{equation}
where $C$ is a positive constant only depending on $|u|_{C^0}$ and other known data, but may change from line to line.

Now we determine the function $h$. Let $\delta$ be a small positive constant to be chosen and $C_0$ be a positive constant such that $ |u|_{C^0} + 1\leqslant C_0$. Let $h$ be defined by
$$
h(u) = e^{\delta(u + C_0)^2}.
$$
Differentiate $h$ twice we have
$$
h^{'} = 2 \delta (u + C_0) e^{\delta(u + C_0)^2} = 2 \delta ( u + C_0) h
$$
and
$$
h^{''} = 4 \delta^2 (u+C_0)^ 2h + 2 \delta h.
$$
It follows that
$$
\frac{h^{''}} {h} - 3 \frac{h^{'2}} {h^2} = 2 \delta - 8 \delta^2 (u + C_0)^2 \geqslant \delta
$$
when $\delta$ is sufficiently small.
By \eqref{45} and \eqref{47}, we can see that
$$
F^{ij} \nabla_{ij}u = F^{ij} U_{ij} - F^{ij} A^{ij} \geqslant F^{ij} A^{ij} \geqslant -C |\nabla_1 u|^ \gamma.
$$
Combining with \eqref{47} and \eqref{32}, we now derive from \eqref{43} that
\begin{equation}
\label{44}
\begin{aligned}
0 \geqslant & \; \delta F^{11} (\nabla_1 u)^2 - C |\nabla_{1} u|^{1 + \gamma_0} (1 + \sum F^{ii}) - C|\nabla_1 u|^ \gamma (1 + \sum F^{ii})\\
& \; - \frac{2Cr} {\zeta} \nabla_{1} u  - \frac{2Cr} {\zeta} |\nabla_{1} u|^{\gamma - 1} \sum F^{ii}-C |\nabla_1 u|^\gamma - \frac{15r^2} {\zeta^2}\sum F^{ii}.
\end{aligned}
\end{equation}
From \eqref{37} we have
\[
\begin{aligned}
U_{11} = &\; - \nabla_{1} u \left( \frac{h^{'}} {h} \nabla_{1} u - \frac{2\rho} {\zeta} \nabla_{1} \rho \right) + A^{11}\\
= &\; -2 \delta (u + C_0) (\nabla_{1} u)^2 + \frac{2\rho} {\zeta} \nabla_{1} u \nabla_{1} \rho + A^{11}.
\end{aligned}
\]
Since $u + C_0 \geqslant 1$, we assume at $x_0$ that $\nabla_{1} u$ is sufficiently large such that $\zeta \nabla_{1} u > \frac{2 \rho \nabla_{1} \rho} {\delta ( u + C_0)}$. This implies that
$U_{11} < - \frac{\delta}{2}|\nabla u|^2$ at $x_0$. By Lemma 1.19 in \cite{Dinew} and the assumption $\tilde f \leqslant C|\nabla u|^\gamma$
in \eqref{47} when $|\nabla u|$ is sufficiently large, we have
$$
F^{11} \geqslant  \nu_0 \sum F^{ii}.
$$
Thus, when $\nabla_{1} u(x_0)$ is sufficiently large, from \eqref{44} we get
\[
\begin{aligned}
0 \geqslant &\; \nu_0 \delta (\nabla_1 u)^2 \sum F^{ii} - C ( |\nabla_1 u|^{1 + \gamma_0} + |\nabla_1 u|^{\gamma} ) \sum F^{ii}\\
%&\; - C \left( 1+ |\nabla_1 u|^{1+\gamma_0} - |\nabla_1 u|^{\gamma} \right) %\sum F^{ii}\\
&\; -\frac{2Cr} {\zeta} \nabla_{1} u - \frac{2Cr} {\zeta} |\nabla_{1} u|^{\gamma - 1} \sum F^{ii} - \frac{15r^2} {\zeta^2} \sum F^{ii},
\end{aligned}
\]
since $\sum F^{ii} \geqslant F(\textbf{1})$.
We therefore obtain that
$$
\zeta(x_0) \nabla_{1} u(x_0) \leqslant \frac{Cr} {\delta \nu_0}.
$$
We just proved that $G(x_0)$ is bounded by some uniform positive constant $C$, which depends only on $|u|_{C^0}$ and other known data. Therefore, by
$ G(0) \leqslant G(x_0)$, we have
$$
|\nabla u(0) |\leqslant \frac{C} {r},
$$
where $C$ depends only on $|u|_{C^0}$ and other known data.
\end{proof}

\begin{remark}
Note that $F(\lambda)$ satisfies \eqref{29},  \eqref{30} and \eqref{45}. Since
$$F_j(\lambda) \geqslant \nu_0 (1 + \sum F_i(\lambda)), \; \mbox{if} \;  \lambda_j<0, \; \forall \; \lambda \in  \mathcal{P}_p$$
 for some positive constant $\nu_0$, the following global gradient estimate
\begin{equation}
\label{49}
\sup_{\bar M} |\nabla u| \leqslant C(1 + \sup_{\partial M} |\nabla u|  )
\end{equation}
can be proved by similar argument as Theorem 5.1 in \cite{GUAN2} without the assumption $|\tilde f| \leqslant C|\nabla u|^\gamma$ in \eqref{47}.
\end{remark}

Let $u \in C^2 (M) \cap C^0 (\bar M)$ be an admissible solution of \eqref{eqn} and \eqref{eqn-b} with $u \geqslant \ul u$, we have
$$\Delta u + \mbox{tr} \chi > 0 \; \mbox{in} \; M.$$
Let $h$ be the solution of
$$\Delta h + \mbox{tr} \chi =0 \; \mbox{in} \; M, \;  h=\varphi \;  \mbox{on}\; \partial M.$$
By the maximum principle,
$$\ul u \leqslant u \leqslant h  \; \mbox{in} \;  M, \; \ul u = u = h \;  \mbox{on} \;\partial M.$$
Consequently, $\ul u_{\nu} \leqslant u_{\nu} \leqslant h_{\nu}$ on $\partial M$, where $\nu$ is the unit inner normal to $\partial M$.
Therefore, we obtain
\begin{equation}
\label{54}
\sup_{\bar M} |u| + \sup_{\partial M} |\nabla u| \leqslant C.
\end{equation}

Combining \eqref{49} and \eqref{54}, we have the $C^1$ estimate as follows.
\begin{theorem}
Assume that \eqref{27} holds for $0 < \gamma_0 < 1$. Let $u \in C^3 (M) \cap C^1 (\bar M)$ be an admissible solution of \eqref{eqn} and \eqref{eqn-b} with $u \geqslant \ul u$ in $M$, where $\ul u$ is a subsolution.
Then, the estimate
\[
	|\nabla u|_{C^0(M)} \leqslant C
\]
 holds for a positive constant $C$ depending on $|u|_{C^0}$, $|\ul u|_{C^1(\bar M)}$  and other known data.
\end{theorem}

\section{Global $C^2$ estimates}
In this section, we establish the global $C^2$ estimate for the equation \eqref{eqn}.
Our main result is the following theorem.
\begin{theorem}
 Assume that there exists an admissible subsolution $\ul u \in C^2 (\bar M)$ satisfying \eqref{26}. Let $u \in C^4 (M) \cap C^2 (\bar M)$ be an admissible solution to \eqref{eqn} on $\bar M$. Then there exists a positive constant $C$ depending on $|u|_{C^1 (\bar M)}$, $|\ul u|_{C^2 (\bar M)}$, $|f|_{C^2}$ and $\inf f$ such that
\begin{equation}
\label{C2}
\sup_{\bar M} |\nabla^2 u| \leqslant C (1 + \sup_{\partial M} | \nabla^2 u|)
\end{equation}
as long as $p \geqslant \frac{n} {2}$.
\end{theorem}
First since
$\lambda (\nabla ^2 \ul u + \chi) \in \mathcal{P}_p$, we can find a positive constant $\epsilon_0$ such that $\lambda (\nabla^2 \ul u + \chi - \epsilon_0 g) \in \mathcal{P}_p$ for all $x \in \bar M$. Note that if $|u|_{C^1 (\bar M)}$ is under control there exist uniform constants
$ \widetilde{f_2} \geqslant \widetilde{f_1} > 0$ such that
\begin{equation}
\label{33}
\begin{aligned}
\widetilde{f_1} \leqslant \widetilde{f} (x, u, \nabla u) \leqslant \widetilde{f_2} \; \mbox{on} \; \bar M.
\end{aligned}
\end{equation}
Hence, by the concavity of $F$, we have
\begin{equation}
\begin{aligned}
\label{5}
F^{ij} \nabla_{ij}(\ul u-u)
\geqslant  F(\ul U - \epsilon_0 g) - F(U) + \epsilon_0 \sum F^{ii}
\geqslant  -C + \epsilon_0 \sum F^{ii},
\end{aligned}
\end{equation}
where $\ul U_{ij} = \nabla_{ij} \ul u + \chi_{ij}.$

Now we are ready to prove Theorem 4.1.
\begin{proof}
Consider the quantity
$$
W= \mbox{max}_{ \xi \in T_x M, |\xi| = 1} \; U_{\xi \xi} e^{\phi},
$$
where $U_{\xi \xi} = U(\xi, \xi)$,
$$
\phi = \frac{\delta |\nabla u|^2} {2} + b(\ul u-u)
$$
and $\delta$ (sufficiently small), $b$ (sufficiently large) are positive constants to be determined.
We may assume $W$ is achieved at an interior point $x_0 \in M$ and $\xi_0 \in T_{x_0} M$. Choose a smooth orthonormal local frame $e_1, \cdots, e_n$  about $x_0$ such that $e_1 = \xi_0$, $\nabla_{e_i} e_j = 0$ and $\{ U_{ij} \}$ is diagonal. So  $\{F^{ij} (x_0)\}$ is diagonal. We  may assume
$$
U_{11} (x_0) \geqslant \cdots \geqslant U_{nn}(x_0).
$$
We have, at $x_0$ where the function $\mbox{log} U_{11} + \phi$ attains its maximum,
\begin{equation}
\label{1}
\begin{aligned}
\frac{\nabla_i U_{11}} {U_{11}} + \delta \nabla_l u \nabla_{il} u + b \nabla_i (\ul u - u) = 0,
\end{aligned}
\end{equation}
for each $i = 1, 2, \cdots, n$ and
\begin{equation}
\label{3}
\begin{aligned}
F^{ii} \left \{ \frac{\nabla_{ii} U_{11}} {U_{11}} - \left( \frac{\nabla_{i} U_{11}} {U_{11}} \right)^2 + \delta (\nabla_{ii} u)^2 + b \nabla_{ii}(\ul u-u) + \delta \nabla_{l} u \nabla_{iil} u \right\} \leqslant 0.
\end{aligned}
\end{equation}
Differentiating the equation \eqref{eqn} twice we have, at $x_0$,
\begin{equation}
\begin{aligned}
\label{34}
F^{ii} \nabla_l {U_{ii}} = \nabla_{l}^{'} \tilde f + \tilde f_u \nabla_l u + \tilde f_{pk} \nabla_{lk}u,
\end{aligned}
\end{equation}
where $\nabla_{l}^{'} \tilde f$ denotes the partial covariant of $\tilde f$ when viewed as depending on $x \in M$ only for $l=i, \ldots, n$. By \eqref{18}, \eqref{1} and  \eqref{34}, we obtain
\begin{equation}
\label{4}
\begin{aligned}
F^{ii} \nabla_{11} {U_{ii}}
\geqslant -F^{ij,kl} \nabla_{1} U_{ij} \nabla_{1} U_{kl} -C U_{11}^2 - C b U_{11}
\end{aligned}
\end{equation}
provided $U_{11}$ is sufficiently large. By \eqref{25},
\begin{equation}
\label{6}
\begin{aligned}
F^{ii} \nabla_{ii} U_{11} - F^{ii} \nabla_{11}{U_{ii}} \geqslant - C U_{11}\sum F^{ii}.
\end{aligned}
\end{equation}
By \eqref{18} and \eqref{34},  we have
\begin{equation}
\label{35}
\begin{aligned}
F^{ii} \nabla_{l} u \nabla_{iil} u
\geqslant &\;- C U_{11} - C \sum F^{ii}.
\end{aligned}
\end{equation}
Now we can derive from \eqref{5}, \eqref{3}, \eqref{4}, \eqref{6} and \eqref{35} that, at $x_0$,
\begin{equation}
\label{7}
\begin{aligned}
0 \geqslant   \frac{\delta} {2} F^{ii} U_{ii}^2 -C U_{11} - F^{ii} \left( \frac{\nabla_{i} U_{11}} {U_{11}} \right)^2
 + \frac{b\epsilon_0} {2}  \sum F^{ii}  - \frac{2} {U_{11}} F^{i1,1i} \nabla_{1} U_{i1} \nabla_{1} U_{1i}
\end{aligned}
\end{equation}
for sufficiently small $\delta $ and $U_{11} \gg b \gg 1$.

We consider two cases: \textbf{(1)}  $U_{nn} \geqslant - \delta_0 U_{11}$ and \textbf{(2)}  $U_{nn} < - \delta_0 U_{11}$, where  $0 < \delta_0 \leqslant \frac{1} {2(p-1)}$.

\textbf{Case (1).}  $U_{nn} \geqslant - \delta_0 U_{11}.$  By Lemma 2.4 in \cite{Dong1}, where used $p\geqslant \frac{n}{2}$, we have
\begin{equation}
\label{8}
\begin{aligned}
\frac{\epsilon_0 b} {4} \sum F^{ii} \geqslant C U_{11}
\end{aligned}
\end{equation}
for sufficiently large $b$. Combining \eqref{7} and \eqref{8}, we have
\[
\begin{aligned}
0 \geqslant  - C b - F^{ii} \left( \frac{\nabla_{i}U_{11}} {U_{11}} \right)^2 + \frac{\delta} {2} F^{ii} U_{ii}^2  + \frac{b\epsilon_0} {4}  \sum F^{ii} - \frac{2} {U_{11}} F^{1i,i1} \nabla_{1} U_{i1} \nabla_{1} U_{1i}.
\end{aligned}
\]
By Lemma 2.5 in \cite{Dong1}, it follows that
\begin{equation}
\label{9}
\begin{aligned}
0\geqslant  - C b - 2F^{11} \left( \frac{\nabla_{1} U_{11}} {U_{11}} \right)^2 + \frac{\delta} {2} F^{ii} U_{ii}^2
+ \frac{b\epsilon_0} {4}  \sum F^{ii}. \\
\end{aligned}
\end{equation}
By \eqref{1}, we obtain
\begin{equation}
\label{52}
\begin{aligned}
\left( \frac{\nabla_{i} U_{11}} {U_{11}} \right)^2
\leqslant &\; C \delta^2 U_{11}^2 + Cb^2.
\end{aligned}
\end{equation}
Substituting \eqref{52} into  \eqref{9}, we have
\begin{equation}
\label{53}
\begin{aligned}
0 \geqslant &\; - C b + \frac{\delta}{4} F^{ii} U_{ii}^2
+ \frac{b\epsilon_0} {4}  \sum F^{ii}
\end{aligned}
\end{equation}
provided  $U_{11}$ is sufficiently large.

\textbf{Case (2).}   $U_{nn} < - \delta_0 U_{11}.$
From \eqref{7}, we have
\begin{equation}
\label{10}
\begin{aligned}
0
%\geqslant &\; -C U_{11} - C b - C \delta^2 F^{ii} U_{ii}^2 - C b^2 \sum %F^{ii} + \frac{\delta} {2} F^{ii} U_{ii}^2 + \frac{b} {2} \epsilon_0 \sum %F^{ii}\\
%\geqslant &\; -C U_{11} + (\frac{\delta} {2}- C \delta^2) F^{ii} U_{ii}^2 - %C b^2 \sum F^{ii} \\
\geqslant &\; -C U_{11} + (\frac{\delta} {2}- C \delta^2) \frac{\delta_0^2 U_{11}^2}{n}  \sum F^{ii} - C b^2 \sum F^{ii}
\end{aligned}
\end{equation}
since $F^{nn} \geqslant \frac{1}{n} \sum F^{ii}$, provided $U_{11} \geqslant  b$ and $\delta$ is sufficiently small.

Now in view of \eqref{53} and \eqref{10} in both cases, we can choose $0 < \delta \ll 1 \ll b$ to obtain
$$
U_{11} \leqslant  C
$$
by the facts $F^{11}U_{11} \geqslant c_0$ and $\sum F^{ii} \geqslant F(\textbf{1})$, respectively.
Therefore, Theorem 4.1 is proved.
\end{proof}

\section{Boundary $C^2$ estimates}
In this section, we establish the boundary $C^2$ estimate for \eqref{eqn}. Our main result is the following theorem.
\begin{theorem}
\label{th}
 Assume that there exists an admissible subsolution $\ul u \in C^2 (\bar M)$ satisfying \eqref{26} and \eqref{eqn-b}. Let $u \in C^3(M) \cap C^2 (\bar M)$ be an admissible solution of \eqref{eqn} and \eqref{eqn-b} with $u \geqslant \ul u $ on $\bar M$. Then, there exists a positive constant $C$  depending on $|u|_{C^1(\bar M)}$, $|\ul u|_{C^2 (\bar M)}$, $|\varphi|_{C^4(\partial M)}$, $|f|_{C^1}$ and $\inf f$ such that
$$
\max_{\partial M} |\nabla^2 u| \leqslant C.
$$
\end{theorem}
For a point $x_0$ on $\partial M$, we choose a smooth orthonormal local frame $e_1, e_2, \dots, e_n$ around $x_0$ such that when restricted to  $\partial M$, $e_n$ is the interior normal to $\partial M$. For $x\in \bar M$ let $\rho(x)$ and $d(x)$ denote the distances from $x$ to $x_0$ and $\partial M$, respectively,
$$
\rho(x):= \mbox{dist}_M(x, x_0), \;\; d(x): = \mbox{dist}_M (x, \partial M).
$$
We may assume $\delta>0$ sufficiently small such that $\rho$ and $d$ are smooth in  $M_{\delta}=\{x\in M : \rho(x)< \delta\}$.
It is easy to estimate the pure tangential second derivatives. Let $e_\alpha, e_\beta$ be some tangential unit vector fields on $\partial M$.
Since $u - \ul u = 0$ on $\partial M$ we have
\begin{equation}
\begin{aligned}
\label{14}
\nabla_{\alpha \beta}(u - \ul u) = -\nabla_n(u-\ul u) \Pi (e_\alpha , e_\beta), \;\; \forall \; 1 \leqslant \alpha, \beta \leqslant n-1 \; \mbox{on}\; \partial M,
\end{aligned}
\end{equation}
where $\Pi$ denotes the second fundamental form of $\partial M$. Therefore,
$$
| \nabla_{\alpha \beta} u (x_0) | \leqslant  C, \;\;  \forall \; 1  \leqslant \alpha, \beta \leqslant n-1.
$$

Next, we prove a bound for the mixed tangential-normal derivatives
\begin{equation}
\label{55}
| \nabla_{\alpha n} u (x_0) | \leqslant C, \;\;  \forall \; 1 \leqslant \alpha \leqslant n-1.
\end{equation}
Define the linear operator $\mathcal{L}$ locally  by
$$
\mathcal{L} w: = F^{ij} \nabla_{ij} w -  \tilde f_{pl} (x, u, \nabla u) \nabla_{l} w, \; w \in C^2(M).
$$
To proceed we first present the key lemma, whose proof is similar to Lemma 3.3 in \cite{D2}.

\begin{lemma}
\label{lem5}
There exist  some uniform positive constants $\varepsilon$, $t$, $\delta$ sufficiently small and $N$ sufficiently large such that the function
$$ v = u - \ul u + t d - Nd^2 $$
satisfies $v \geqslant 0 \; \mbox{on} \; \bar M_\delta$ and
\begin{equation}
\label{50}
\mathcal{L} v \leqslant \; -\varepsilon (1 + \sum F^{ii}) \; \mbox{in} \; M_\delta.
\end{equation}
\end{lemma}
\begin{proof}
By direct calculations, we see that
$$\mathcal{L} v= \mathcal{L} (u - \ul u) + (t-2Nd) \mathcal{L} d - 2NF^{ij} \nabla_i d  \nabla_j d.$$
By the definition of the operator $\mathcal{L}$ and \eqref{5}, we have
\[
\begin{aligned}
\mathcal{L} (u - \ul u)
%= &\; F^{ij} \nabla_{ij} (u - \ul u) - \tilde f_{pl} \nabla_{l} (u - \ul u) \\
\leqslant &\;  C - \varepsilon_0 \sum F^{ii},
\end{aligned}
\]
where $C$ depends on $\sup_M  \widetilde{f}$ and $ \ul u$.
It also easy to see that
$$|\mathcal{L}d| \leqslant  C ( 1+ \sum F^{ii} )$$
for some positive constant $C$ under control. Therefore, we obtain that
$$\mathcal{L}v \leqslant  C - \varepsilon_0 \sum F^{ii}  - 2NF^{ij} \nabla_i d \nabla_j d + C (t + 2Nd) \sum F^{ii}.$$
Choosing $t \ll 1$ and $\delta \ll 1$, we have
$$\mathcal{L} v \leqslant  C - \frac{3\varepsilon_0}{4} \sum F^{ii} - 2 N F^{nn},$$
since $|\nabla_n d|=1$.
Without loss of generality let us assume $F_1 \leqslant \cdots \leqslant F_n$. We have $ F^{nn} \geqslant F_1$.  Therefore, we have
\begin{equation}
\begin{aligned}
\label{19}
4 N F^{nn} + \varepsilon_0 \sum F^{ii}  \geqslant  (4N)^{p/n} \varepsilon_0^{(n-p)/n}
\end{aligned}
\end{equation}
by the inequality of arithmetic and geometric means.
Choosing $N$ sufficiently large, we derive from \eqref{19} that
$$\mathcal{L} v \leqslant - \frac{\varepsilon_0}{4}  \sum F^{ii} \leqslant - \frac{\varepsilon_0}{8} ( 1+ \sum F^{ii} )$$
as $\sum F^{ii} \geqslant p \geqslant 1$.
\end{proof}

With the above lemma in hand, we can construct a suitable barrier function to prove \eqref{55}.
Given $1 \leqslant \alpha \leqslant  n-1$, we define
\begin{equation}
\begin{aligned}
\label{23}
w = \pm\nabla_ \alpha (u - \varphi) - \sum_{l < n} |\nabla_l (u - \varphi)|^2,
\end{aligned}
\end{equation}
where $\varphi$ is extended to $M$ with  $\nabla_n \varphi = 0$ on $\partial M$.
By a straightforward calculation using \eqref{18} and \eqref{34} for each $1 \leqslant k \leqslant n$,
\begin{equation}
\begin{aligned}
\label{11}
|\mathcal{L}\nabla_k (u - \varphi)|
\leqslant &\; C \big( 1 + \sum F^{ii} + C \sum \frac{\partial F} {\partial \lambda_i}| \lambda_i|\big)
\end{aligned}
\end{equation}
and
\begin{equation}
\begin{aligned}
\label{56}
\mathcal{L} (| \nabla_k (u - \varphi) |^2)  \geqslant F^{ij} U_{ik} U_{jk} -  C \big( 1 + \sum F^{ii} +  \sum \frac{\partial F} {\partial \lambda_i}| \lambda_i|\big).
\end{aligned}
\end{equation}
By proposition 2.19 in \cite{G} there exists an index $r$ such that
\begin{equation}
\begin{aligned}
\label{13}
 \sum_{l<n} F^{ij} U_{il} U_{jl} \geqslant \frac{1} {2} \sum_{i\neq r} \frac{\partial F} {\partial \lambda_i} \lambda_i^2.
\end{aligned}
\end{equation}
Thus, by \eqref{56}, we have
\begin{equation}
\begin{aligned}
\label{57}
\sum_{l<n}\mathcal{L} (| \nabla_l (u - \varphi) |^2)  \geqslant \frac{1} {2} \sum_{i\neq r} \frac{\partial F} {\partial \lambda_i} \lambda_i^2 - C \big( 1 + \sum F^{ii} +  \sum \frac{\partial F} {\partial \lambda_i}| \lambda_i|\big).
\end{aligned}
\end{equation}
Combining \eqref{23}, \eqref{11}  and \eqref{57}, we obtain
\[
\begin{aligned}
\mathcal{L} w &\; \leqslant C (1 + \sum F^{ii}) + C \sum \frac{\partial F} {\partial  \lambda_i} |\lambda_i| - \frac{1}{2} \sum_{i\neq r } \frac{\partial F} {\partial \lambda_i} \lambda_i^2
\end{aligned}
\]
for some index $r$. By Corollary 2.21 in \cite{G} or see \cite{CP},
we have,
for any index $1\leqslant r \leqslant n$,
\begin{equation}
	\label{60}
\sum F_i |\lambda_i| \leqslant \epsilon \sum_{i\neq r} F_i \lambda_i^2
+ \frac{C}{\epsilon} \sum F_i + C
\end{equation}
for any $\epsilon >0$.
It follows that
\begin{equation}
\begin{aligned}
\label{17}
\mathcal{L} w &\; \leqslant C (1 + \sum F^{ii}).
\end{aligned}
\end{equation}

Now we fix $\delta, t$ such that \eqref{50} holds. Let $h = w + B \rho^2 + Av$ and choose $A \gg B\gg 1$ such that $h \geqslant 0$ on $\partial M_\delta$,  from  Lemma \ref{lem5} and \eqref{17} we obtain
\[
\begin{aligned}
\mathcal{L} h &\; \leqslant (C + CB - A\varepsilon) (1 + \sum F^{ii}) \leqslant 0 \; \; \mbox{in} \; M_\delta.
\end{aligned}
\]
By the Maximum Principle, we derive  $h \geqslant 0$  in $M_\delta$ and $\nabla_n h(x_0) \geqslant 0$. So we have \eqref{55}.

It remains to derive the double normal estimate
$$
|\nabla_{n n} u (x_0)| \leqslant C.
$$
For $x_0 \in \partial M$, denote $\lambda (x_0) = \lambda (U(x_0)) = (\lambda_1,\lambda_2,\cdots, \lambda_n)$.
Choose a smooth orthonormal local frame $e_1, \cdots, e_n$ around $x_0$ such that $\{U_{\alpha\beta} (x_0)\}$ is diagonal,
%\[
%U(x_0)=\left( \begin{array}{ccccc}
%U_{11}(x_0) & 0 & \cdots & 0& U_{1n}(x_0) \\
%0 & U_{22}(x_0) & \cdots & 0& U_{2n}(x_0) \\
%\vdots & \vdots & \ddots & \vdots & \vdots\\
%0& 0 & \cdots & U_{n-1 n-1}(x_0) & U_{n-1 n}(x_0) \\
%U_{n1}(x_0) & U_{n2}(x_0) & \cdots & U_{n n-1}(x_0) & U_{nn}(x_0) %\end{array} \right)
%\]
%\[
%U(x_0)=\left( \begin{array}{cc}
%	U_{\alpha\beta}(x_0) & U_{\alpha n}(x_0) \\
%	U_{n\beta}(x_0)  & U_{nn}(x_0) \end{array} \right)
%\]
where $1\leqslant \alpha, \beta \leqslant n-1$.
Denote the upper left $(n-1)\times(n-1)$ diagonal matrix in $U(x_0)$ by $\tilde U (x_0)$
%\[\tilde U (x_0) = \left( \begin{array}{ccccc}
%U_{11}(x_0) & \cdots & 0  \\
%\vdots  & \ddots & \vdots   \\
%0& \cdots &  U_{n-1n-1}(x_0) \\
%\end{array} \right)
%\]
and let
\[\tilde \lambda (x_0) = \lambda (\tilde U(x_0)) = (U_{11} (x_0), \cdots, U_{n-1 n-1} (x_0)).\]
Without loss of generality, we assume
$$\lambda_1 \leqslant  \lambda_2 \leqslant  \cdots \leqslant \lambda_n$$
and
$$U_{11}(x_0) \leqslant U_{22}(x_0) \leqslant \cdots \leqslant U_{n-1 n-1}(x_0).$$
By Cauchy interlace theorem, we see
$$\lambda_1 \leqslant U_{11}(x_0) \leqslant \lambda_2 \leqslant \cdots \leqslant \lambda_{n-1} \leqslant U_{n-1n-1}(x_0) \leqslant \lambda_n.$$

For every $x\in \partial M$ with a local frame chosen around $x$ as above,
we shall prove that
\[
U_{11}(x) + U_{22}(x) + \cdots + U_{pp}(x) \geqslant c_0
\]
for some uniform positive constants $c_0$.
We use an idea from Trudinger \cite{Tru}.
It suffices to show that at the point $y_0 \in \partial M$ where the above quantity attains a minimum,
\begin{equation}
\begin{aligned}
\label{51}
U_{11}(y_0) + U_{22}(y_0) + \cdots + U_{pp}(y_0) \geqslant c_0.
\end{aligned}
\end{equation}
Choose a local orthonormal frame $e_1, \cdots, e_n$ about $y_0$ as before such that $\tilde U(y_0)$ is  diagonal.
%Since $\nabla_n \varphi = 0$, we have
%\[\begin{aligned}
%U_{ii}(x) %= &\;\nabla_{ii}u(x) + \chi^{ii} \\
%= &\; \nabla_{ii}\varphi + \nabla_{n}(\varphi-u) \Pi (e_ie_i) + \chi^{ii}\\
%= &\; \nabla_{ii} \varphi - \nabla_{n} u \; \Pi (e_i,e_i) + \chi^{ii}.
%\end{aligned}
%\]
For $x \in \partial M$ near $y_0$, we have
\begin{equation}
\begin{aligned}
\label{15}
\sum_{i=1}^p U_{ii}(x) \geqslant \sum_{i=1}^p U_{ii}(y_0).
\end{aligned}
\end{equation}
By \eqref{14} and \eqref{15}, we can see that  for $x\in \partial M$ near $y_0$,
\begin{equation}
\begin{aligned}
\label{16}
 \nabla_n (u - \ul u) (x) \sum_{i=1}^p \Pi (e_i,e_i)
%= &\; \sum_{i=1}^p \nabla_{ii} (\ul u- u) (x) \\
%= &\; \sum_{i=1}^p \ul U_{ii} (x) - \sum_{i=1}^p U_{ii} (x) \\
\leqslant &\; \sum_{i=1}^p \ul U_{ii} (x) - \sum_{i=1}^p  U_{ii} (y_0),
\end{aligned}
\end{equation}
where $\ul U_{ij} = \nabla_{ij} \ul u + \chi_{ij}.$
Since $\ul u$ is an admissible subsolution and $\lambda(\ul U)\in \mathcal{P}_p$, there exists a uniform constant $c_0{'}>0$ such that
$$
\min_{x\in \partial M} \sum_{i=1}^p \ul U_{ii} (x) \geqslant c_0^{'}.
$$
We may assume $\sum_{i=1}^p U_{ii} (y_0) < \frac{c_0^{'}}{2}$ for otherwise we are done. By \eqref{14}, we have
\[
\begin{aligned}
\nabla_n (u-\ul u) (y_0) \sum_{i=1}^p \Pi (e_i,e_i) = &\; \sum_{i=1}^p \ul U_{ii} (y_0) - \sum_{i=1}^p U_{ii} (y_0) \geqslant \frac{c_0{'}}{2},
\end{aligned}
\]
which implies that
\[
\begin{aligned}
\sum_{i=1}^p \Pi (e_i,e_i)\geqslant c_1 \; \mbox{in} \; \bar M_{\delta^{'}} \\
\end{aligned}
\]
for some uniform positive constant $c_1$ and $\delta^{'} \leqslant \delta$.
%Let $\psi (x)=\sum_{i=1}^p U_{ii} (x)$
%and
Define
\[\begin{aligned}
\Phi(x) = \frac{ \sum_{i=1}^p \ul U_{ii} (x) - \sum_{i=1}^p U_{ii} (y_0) }{\sum_{i=1}^p \Pi (e_i,e_i)},
\end{aligned}
\]
which is smooth in $\bar M_{\delta^{'}}$. Consequently, by \eqref{16} we have
\[
\begin{aligned}
\nabla_n (u - \ul u)(x) \leqslant \Phi(x) \; \mbox{on} \; \partial M \cap \bar M_{\delta^{'}}.
\end{aligned}
\]
Note that $\nabla_n (u - \ul u)(y_0) = \Phi(y_0)$.
Let
\[\begin{aligned}
\Psi = A_1 v +  A_2 \rho^2 - A_3 \sum_{\beta < n} |\nabla_\beta (u - \varphi)|^2.
\end{aligned}
\]
By Lemma 5.2, \eqref{11}, \eqref{56} and \eqref{13}, we have
\[
\begin{aligned}
\mathcal{L} \Psi
%= &\; A_1 F^{ij} \nabla_{ij} v + A_2 2\rho F^{ij} \nabla_{ij} \rho + 2 A_2 %F^{ij} \nabla_{i} \rho \nabla_{j} \rho \\
%&\; - A_3 2 \sum_{\beta<n} |\nabla_\beta (u - \varphi)| F^{ij} \nabla_{ij} %\nabla_\beta (u - \varphi) \\
%&\; - 2 A_3 \sum_{\beta<n} F^{ij} \nabla_{i} \nabla_\beta (u - \varphi) %\nabla_{j} \nabla_ \beta (u - \varphi)\\
%\leqslant &\; -A_1 \varepsilon (1 + \sum F^{ii}) + C A_2 \sum F^{ii} + C %A_3+ C A_3 \sum F^{ii}  \\
%&\; + C A_3 \sum \frac{\partial F} {\partial \lambda_i } |\lambda_i| - A_3 %\sum_{\beta<n} F^{ij} U_{i\beta} U_{j\beta} \\
\leqslant &\; -A_1 \varepsilon (1 + \sum F^{ii}) + C (A_2 + A_3) \sum F^{ii}\\
&\; + C A_3 + C A_3 \sum \frac{\partial F} {\partial \lambda_i} |\lambda_i| - \frac{A_3}{2} \sum_{i\neq r} \frac{\partial F} {\partial \lambda_i } \lambda_i ^2
\end{aligned}
\]
for some index $r$.
By \eqref{60}, it follows that
\begin{equation}
\begin{aligned}
\label{24}
\mathcal{L} \Psi
%\leqslant &\; - A_1 \varepsilon (1 + \sum F^{ii}) + C (A_2 + A_3) \sum %F^{ii} + C A_3 \\
\leqslant - \frac{A_1}{2} \varepsilon (1 + \sum F^{ii})
\end{aligned}
\end{equation}
for  $A_1 \gg A_2 \gg A_3 \geqslant 1$.  Thus, by \eqref{11} and \eqref{24}, we can choose $\Psi$ such that
\[
\begin{aligned}
\mathcal{L} ( \nabla_{n} (u - \ul u) - \Phi - \Psi ) \geqslant 0 \; \mbox{in} \; M_{\delta_0}, \\
\end{aligned}
\]
and
\[
\begin{aligned}
\nabla_{n} (u - \ul u)- \Phi - \Psi \leqslant 0 \; \mbox{on} \; \partial M_{\delta_0}
\end{aligned}
\]
for $A_1 \gg A_2 \gg A_3 \geqslant 1$ and some  positive constant $\delta_0 \leqslant \delta^{'}.$ By the maximum principle, we get a bound
\[
\begin{aligned}
\nabla_{nn} u(y_0) \leqslant  C. \\
\end{aligned}
\]
Then, at $y_0$, $\lambda (U) (y_0) \in \mathcal{P}_p$ is bounded. Therefore,
\[
\begin{aligned}
\sum_{i=1}^p \lambda_{i} (y_0) \geqslant c_0  \\
\end{aligned}
\]
for some $c_0>0$. This implies \eqref{51} holds.

Now we can prove the estimates for the double normal derivatives.
By Lemma 1.2 in \cite{CNS3}, we have
$$
\lambda_1 + \lambda_2 + \cdots + \lambda_p > \frac{c_0}{2}
$$
if $ U_{nn} (x_0) $ is large enough.
Thus,
\[
\begin{aligned}
f= &\; \Pi_{1 \leqslant i_1 < \cdots < i_{p} \leqslant n-1} (\lambda_{i_1} + \cdots + \lambda_{i_{p}})\\
\geqslant &\; (\lambda_1 + \lambda_2 + \cdots +\lambda_p)^{ {C_{n}^{p}} - {C_{n-1}^{p-1}} } \Pi_{1 \leqslant i_1 < \cdots < i_{p} \leqslant n-1}  (\lambda_{n} + \lambda_{i_1} + \cdots + \lambda_{i_{p-1}}) \\
\geqslant &\; (\frac{c_0}{2})^{{C_{n}^{p}} - {C_{n-1}^{p-1}}} \Pi_{1\leqslant i_1 < \cdots < i_{p} \leqslant n-1} (\lambda_{n} + \lambda_{i_1} + \cdots + \lambda_{i_{p-1}} ),
\end{aligned}
\]
which indicates
$$\lambda_n \leqslant  C.$$
Therefore, we have
$$
U_{nn}(x_0) \leqslant C.
$$
The proof of Theorem 5.1 is complete.

\textbf{Data Availability}
Data sharing is not applicable to this article as no datasets were generated or analyzed during the current study.

\textbf{Conflict of Interest}
The authors declare that they have no conflict of interest.

\end{document}